\newcommand{\dsp}{\renewcommand{\baselinestretch}{1.05}}

\documentclass[pdflatex,fleqn]{article}
\usepackage{}
\usepackage{xcolor}
\usepackage{mathrsfs}
\usepackage{amsfonts}
\usepackage{amssymb}
\usepackage{latexsym}
\usepackage{amsmath}
\usepackage{amscd}
\usepackage{graphicx}
\usepackage{hyperref}
\usepackage{hypernat}
\usepackage{cite}

\usepackage{geometry,amsthm,graphics,amssymb,amsmath,enumerate,latexsym,tabularx,shapepar}
\usepackage[all,2cell,dvips]{xy}\UseAllTwocells\SilentMatrices

\usepackage{extarrows}
\date{}

\dsp

\title{\bf Reduction to dimension two of local spectrum for $AH$ algebra with ideal property}

\author{
 Chunlan Jiang
}

\begin{document}

\large

\maketitle

\noindent{\underline {\noindent{ }\hspace {157mm}}}\\
Abstract: A $C^{*}$-algebra $A$ has ideal property if any ideal $I$ of $A$ is generated as a closed two sided ideal by the projections inside the ideal. Suppose that the limit $C^{*}$-algebra $A$ of inductive limit of direct sums of matrix algebras over spaces with uniformly bounded dimension has ideal property. In this paper we will prove that $A$ can be written as an inductive limit of certain very special subhomogeneous algebras, namely, direct sum of dimension drop interval algebras and matrix algebras over 2-dimensional spaces with torsion $H^{2}$ groups. 


\vspace{3mm}

\noindent\textbf{\S1. Introduction}

\vspace{-2.4mm}

An $AH$ algebra is a nuclear $C^{*}$-algebra of the form $A=\lim\limits_{\longrightarrow}(A_{n}, \phi_{n,m})$ with $A_{n}=\bigoplus_{i=1}^{t_{n}}P_{n,i}M_{[n,i]}(C(X_{n,i}))P_{n,i}$, where $X_{n,i}$ are compact metric spaces, $t_{n}$,  $[n,i]$ are positive integers, $M_{[n,i]}(C(X_{n,i}))$ are algebras of $[n,i]\times[n,i]$ matrices with entries in $C(X_{n,i})$---the algebra of complex - valued functions on $X_{n,i}$, and finally $P_{n,i}\in M_{[n,i]}(C(X_{n,i}))$ are projections (see [Bla]). If we further assume that $\sup\limits_{n,i}\dim(X_{n,i})<+\infty$ and $A$ has ideal property---each ideal $I$ of $A$ is generated by the projections inside the ideal, then it is proved in [GJLP1-2] that $A$ can be written as inductive limit of $B_{n}=\bigoplus_{i=1}^{s_{n}}P^{'}_{n,i}M_{[n,i]^{'}}(C(Y_{n,i}))P^{'}_{n,i}$. 

In this paper, we will further reduce the dimension of local spectra (that is the spectra of $A_{n}$ or $B_{n}$ above) to 2 (instead of 3). Namely, the above $A$ can be written as inductive limit of direct sum of matrix algebras over the $\{pt\}, [0,1], S^{1}, T_{\uppercase\expandafter{\romannumeral2}, k}$ (no $T_{\uppercase\expandafter{\romannumeral3}, k}$ and $S^{2}$) and $M_{l}(I_{k})$, where $I_{k}$ is dimension drop interval algebra,
$$I_{k}=\{f\in C([0,1],M_{k}(\mathbb{C})), f(0)=\lambda\textbf{1}_{k}, f(1)=\mu\textbf{1}_{k}, \lambda,\mu\in\mathbb{C}\}.$$
In this paper, we will also call $\bigoplus^{s}_{i=1}M_{l_{i}}(I_{k_{i}})$ a dimension drop-algebra.\\
This result unifies the theorems of [DG] and [EGS] (for the rank zero case) and [Li4] (for the simple case). Note that Li$^{'}$s reduction theorem was not used in the classification of simple $AH$ algebra, and Li$^{'}$s proof depends on the classification of simple $AH$ algebra (see [Li4] and [EGL]). For our case, the reduction theorem is an important step toward the classification (see [GJL]). The proof is more difficult than Li$^{'}$s case. For example, in the case of $AH$ algebra with ideal property, one can not remove the space $S^{2}$ without introduce $M_{l}(I_{k})$ (for the simple case, the space $S^{2}$ is removed from the list of spaces in [EGL] without introduced dimension drop algebras). Another point is that, in the simple $AH$ algebras, one can assume each partial map $\phi^{i,j}_{n,m}$ is injective, but in $AH$ algebras with ideal property, we can not make such assumption. For the classification of real rank zero $AH$ algebras, we refer to the readers [Ell1], [EG1-2], [D1-2], [G1-4] and [DG]. For the classification simple $AH$ algebra, we refer to the readers [Ell2-3], [Li1-3], [G5] and [EGL1-2].

The paper is organized as follows.\\
In section 2, we will do some necessary preparation. In section 3, we will prove our main theorem.

\vspace{3mm}

\noindent\textbf{\S2. Preparation}

\vspace{-2.4mm}

We will adopt all the notation from section 2 of [GJLP2]. For example we refer the reader to [GJLP2] for the concepts of $G-\delta$ multiplicative maps (see definition 2.2 there), spectral variation $SPV(\phi)$ of a homomorphism $\phi$ (see 2.12 there) weak variation $\omega(F)$ of a finite set $F\subset QM_{N}(C(X))Q$ (see 2.16 there).\\
As in 2.17 of [GJLP2], we shall use $\bullet$ to denote any possible integer.

\noindent\textbf{2.1.} In this article, without lose of generality we will assume the $AH$ algebras $A$ are inductive limit of
$$A=\lim(A_{n}=\bigoplus_{i=1}^{t_{n}}M_{[n,i]}(C(X_{n,i})), \phi_{n,m}),$$
where $X_{n,i}$ are the spaces of $\{pt\}, [0,1], S^{1}, T_{II, k}, T_{III, k}$, and $S^{2}$. (Note by the main theorem of [GJLP2], all $AH$ algebras with ideal property and with no dimension growth are corner subalgebras of the above form (also see 2.7 of [GJLP2]).)

\noindent\textbf{2.2.} Recall a projection $P\in M_{k}(C(X))$ is called a trivial projection if it is unitarily\\
equivalent to  $\begin{pmatrix}\textbf{1}_{k_{1}}& 0\\0 & 0\end{pmatrix}$ for $k_{1}=rank(P)$. If $P$ is a trival projection and $rank(P)=k_{1}$, then
$$PM_{k}(C(X))P\cong M_{k_{1}}(C(X)).$$

\noindent\textbf{2.3.} Let $X$ be a connected finite simplicial complex, $A=M_{k}(C(X))$, a unital $\ast$ homomorphism $\phi: A\longrightarrow M_{l}(A)$ is called a (unital) simple embedding if it is homotopic to the homomorphism $id\oplus\lambda$, when $\lambda: A\longrightarrow M_{l-1}(A)$ is defined by
$$\lambda(f)=diag(\underbrace{f(x_{0}),f(x_{0}),\cdots,f(x_{0})}\limits_{l-1}),$$
for a fixed base point $x_{0}\in X$.

The following two lemmas are special cases of Lemma 2.15 of [EGS] (also see 2.12 of [EGS]).

\noindent\textbf{Lemma 2.4.} (c.f 2.12 or case 2 of 2.15 in [EGS]). For any finite set $F\subset A=M_{n}(C(T_{\uppercase\expandafter{\romannumeral3}, k}))$ and $\varepsilon>0$, there is a unital simple embedding $\phi: A\longrightarrow M_{l}(A)$ (for $l$ large enough) and a $C^{*}$-algebra $B\subset A$, which is a direct sum of dimension drop algebras and a finite dimensional $C^{*}$-algebra such that
$$dist(\phi(f),B)<\varepsilon,~~~~~~\forall\;f\in F.$$

\noindent\textbf{Lemma 2.5.}(see case 1 of 2.15 in [EGS]) For any finite set $F\subset M_{n}(C(S^{2}))$ and $\varepsilon>0$, there is a unital simple embedding $\phi:A\longrightarrow M_{l}(A)$ (for $l$ large enough) and a $C^{*}$-algebra $B\subset A$, which is a finite dimensional $C^{*}$-algebra such that
$$dist(\phi(f),B)<\varepsilon,~~~~~~\forall\;f\in F.$$

The following lemma is well known:

\noindent\textbf{Lemma 2.6.}(see [G5,4.40]) For any $C^{*}$-algebra $A$ and finite set $F\subset A$, $\varepsilon>0$, there is a finite set $G\subset A$ and $\eta>0$ such that if $\phi: A\longrightarrow B$ is a homomorphism and $\psi: A\longrightarrow B$ is a completely positive linear map, satisfing the following
$$\|\phi(g)-\psi(g)\|<\eta,~~~~~\forall\;g\in G,$$
then $\psi$ is the $F-\varepsilon$ multiplicative.\\

\noindent\textbf{Lemma 2.7.} Let $A=M_{n}(C(T_{\uppercase\expandafter{\romannumeral3}, k}))$ or $M_{n}(C(S^{2}))$. And let a finite set $F\subset A$ and $\varepsilon>0$, there is a diagram
$$
\xymatrix{
A\ar[rr]^{\phi}\ar[rdrd]^{\beta}  &  &  M_{l}(A)  &  &  \\
     &  &  &  &  &  &  \\
                                                                     &  &   B\ar[uu]_{\iota} & &   \\
 }
$$
with the following conditions.\\
(1) $\phi$ is a simple embedding,\\
(2) if $A=M_{n}(C(S^{2}))$, then $B$ is a finite dimensional $C^{*}$-algebra, and if $A=M_{n}(C(T_{\uppercase\expandafter{\romannumeral3}, k}))$, then $B$ is a direct sum of dimension drop $C^{*}$-algebras and a finite dimensional $C^{*}$-algebra, and $\iota$ is an inclusion,\\
(3) $\|\iota\circ\beta(f)-\phi(f)\|<\varepsilon$, $\forall f\in F$, and $\beta$ is $F-\varepsilon$ multiplicative.
\begin{proof}
 Let $G$ and $\eta$ be as Lemma 2.6 for $F$ and $\varepsilon$. Apply Lemma 2.4 or Lemma 2.5 to $A$, $F\cup G\subset A$ and $\frac{1}{3}\min(\varepsilon,\eta)$. One can find a unital simple embedding $\phi: A\longrightarrow M_{l}(A)$, and an sub-$C^{*}$-algebra $B\subset M_{l}(A)$ as required in condition (2) such that
$$dist(\phi(f),B)<\frac{1}{3}\min(\varepsilon,\eta),~~~~~~\forall\;f\in F.$$
Choose a finite $\widetilde{F}\subset B$ such that
$$dist(\phi(f),\widetilde{F})<\frac{1}{3}\min(\varepsilon,\eta),~~~~~~\forall\;f\in F.$$
Since $B$ is a nuclear $C^{*}$-algebra, there are two completely positive linear maps
$$\lambda_{1}: B\longrightarrow M_{N}(\mathbb{C})~~~\mbox{and}~~~\lambda_{2}: M_{N}(\mathbb{C})\longrightarrow B,$$
such that
$$\|\lambda_{2}\circ\lambda_{1}(g)-g\|<\frac{1}{3}\min(\varepsilon,\eta),~~\forall g\in\widetilde{F}$$
Using Arveson extension theorem, one can extend $\lambda_{1}: B\longrightarrow M_{N}(\mathbb{C})$ to a map\\
 $\beta_{1}: M_{l}(A)\longrightarrow M_{N}(\mathbb{C})$. Then it is a straight forward to prove that
$$\beta=\lambda_{2}\circ\beta_{1}\circ\phi: A\longrightarrow B,$$
is as desired.
\end{proof}

The following is modification of Theorem 3.8 of [GJLP2].

\noindent\textbf{Proposition 2.8.} Let $\lim\limits_{n\rightarrow\infty}(A_{n}=\bigoplus_{i=1}^{t_{n}}M_{[n,i]}(C(X_{n,i})), \phi_{n,m})$ be $AH$ inductive limit with ideal property, with $X_{n,i}$ being one of $\{pt\}, [0,1], S^{1}, T_{II, k}, T_{III, k}$, or $S^{2}$. Let $B=\bigoplus_{i=1}^{s}B^{i}$, where $B^{i}=M_{l_{i}}(C(Y_{i}))$, with $Y_{i}$ being  $\{pt\}, [0,1], S^{1}$, or $T_{II, k}$, (no $T_{\uppercase\expandafter{\romannumeral3}, k}$ or $S^{2}$) or $B^{i}=M_{l_{i}}(I_{k_{i}})$---dimension drop $C^{*}$-algebras. Suppose that
$$\widetilde{G}(=\oplus\widetilde{G}^{i})\subset G(=\oplus G^{i})\subset B(=\oplus B^{i}),$$
is a finite set, $\varepsilon_{1}$ is a positive number with $\omega(\widetilde{G}^{i})<\varepsilon_{1}$, if $Y_{i}=T_{\uppercase\expandafter{\romannumeral2}, k}$, and $L$ is any positive integer. Let $\alpha: B\longrightarrow A_{n}$ be any homomorphism. Denote
$$\alpha(\textbf{1}_{B}):= R(=\oplus R^{i})\in A_{n}(=\oplus A^{i}_{n}).$$
Let $F\subset RA_{n}R$ be any finite set and $\varepsilon<\varepsilon_{1}$ be any positive number. It follows that there are $A_{m}$, and mutually orthogonal projections $Q_{0},Q_{1},Q_{2}\in A_{m}$ with
$$\phi_{n,m}(R)=Q_{0}+Q_{1}+Q_{2},$$
a unital map $\theta_{0}\in Map(RA_{n}R,Q_{0}A_{m}Q_{0})_{1}$, two unital homomorphisms $\theta_{1}\in Hom(RA_{n}R,Q_{1}A_{m}Q_{1})_{1}$, $\xi\in Hom(RA_{n}R,Q_{2}A_{m}Q_{2})_{1}$ such that \\
(1) $\parallel \phi_{n,m}(f)-(\theta_{0}(f)\oplus\theta_{1}(f)\oplus\xi(f))\parallel<\varepsilon$, for all $f\in F$,\\
(2) there is a unital homomorphism
$$\alpha_{1}: B\longrightarrow(Q_{0}+Q_{1})A_{m}(Q_{0}+Q_{1}),$$
such that
$$\parallel \alpha_{1}(g)-(\theta_{0}+ \theta_{1})\circ\alpha(g)\parallel<3\varepsilon_{1}~~~~\forall g\in \widetilde{G}_{i},~~~~~~~\mbox{if}~~ B^{i} ~~\mbox{is of form }~~M_{\bullet}(T_{\uppercase\expandafter{\romannumeral2}, k})$$
and
$$\parallel \alpha_{1}(g)-(\theta_{0}+\theta_{1})\circ\alpha(g)\parallel<\varepsilon,~~~~\forall g\in G^{i},~~~~~~~\mbox{if}~~ B^{i} ~~\mbox{is not of form }~~M_{\bullet}(T_{\uppercase\expandafter{\romannumeral2}, k}).$$
(3) $\theta_{0}$ is $F-\varepsilon$ multiplicative and $\theta_{1}$ satisfies that
$$  \theta^{i,j}_{1}([e])\geqslant L\cdot[\theta^{i,j}_{0}(R^{i})].$$
(4) $\xi$ factors through a $C^{\ast}$-algebra C---a direct sum of matrix algebras over C[0,1] as
$$\xi:RA_{n}R\xrightarrow{\xi_{1}}C\xrightarrow{\xi_{2}}Q_{2}A_{m}Q_{2}.$$
\noindent\textbf{Proposition 2.9.} Let $\lim\limits_{n\rightarrow\infty}(A_{n}=\bigoplus_{i=1}^{t_{n}}M_{[n,i]}(C(X_{n,i})), \phi_{n,m})$ be $AH$ inductive limit with ideal property, with $X_{n,i}$ being one of $\{pt\}, [0,1], S^{1}, T_{\uppercase\expandafter{\romannumeral2}, k}, T_{\uppercase\expandafter{\romannumeral3}, k}$ or $S^{2}$. Let $B=\bigoplus_{i=1}^{s}B^{i}$, where $B^{i}=M_{l_{i}}(C(Y_{i}))$, with $Y_{i}$ being  $\{pt\}, [0,1], S^{1}$ or $T_{\uppercase\expandafter{\romannumeral2}, k}$, (no $T_{\uppercase\expandafter{\romannumeral3}, k}$ or $S^{2}$) or $B^{i}=M_{l_{i}}(I_{k_{i}})$---dimension drop $C^{*}$-algebras. Suppose that
$$\widetilde{G}(=\oplus\widetilde{G}^{i})\subset G(=\oplus G^{i})\subset B(=\oplus B^{i}),$$
is a finite set, $\varepsilon_{1}$ is a positive number with $\omega(\widetilde{G}^{i})<\varepsilon_{1}$, if $Y_{i}=T_{\uppercase\expandafter{\romannumeral2}, k}$, and $L>0$ is any positive integer. Let $\alpha: B\longrightarrow A_{n}$ be any homomorphism. Let $F\subset A_{n}$ be any finite set and $\varepsilon<\varepsilon_{1}$ be any positive number. It follows that there are $A_{m}$ and mutually orthogonal projections $P,Q\in A_{m}$ with $\phi_{n,m}(\textbf{1}_{A_{n}})=P+Q$,
a unital map $\theta\in Map(A_{n},PA_{m}P)_{1}$, and a unital homomorphism $\xi\in Hom(A_{n},QA_{m}Q)_{1}$ such that \\
(1) $\parallel \phi_{n,m}(f)-(\theta(f)\oplus\xi(f))\parallel<\varepsilon$, for all $f\in F$,\\
(2) there is a homomorphism $\alpha_{1}: B\longrightarrow PA_{m}P$ such that
\begin{center}
$\parallel \alpha^{i,j}_{1}(g)-(\theta\circ\alpha)^{i,j}(g)\parallel<3\varepsilon_{1}$~~~ $\forall g\in \widetilde{G}^{i}$,
~~~~~if $B^{i}$~~ is of form~~ $M_{\bullet}(C(T_{\uppercase\expandafter{\romannumeral2}, k}))$,
\end{center}
\begin{center}
$\parallel \alpha^{i,j}_{1}(g)-(\theta\circ\alpha)^{i,j}(g)\parallel<\varepsilon$~~~ $\forall g\in G^{i}$, ~~~~~if $B^{i}$~~ is not of form ~~$M_{\bullet}(C(T_{\uppercase\expandafter{\romannumeral2}, k})).$
\end{center}
(3) $\omega(\theta(F))<\varepsilon$ and $\theta$ is $F-\varepsilon$ multiplicative.\\
(4) $\xi$ factors through a $C^{\ast}$-algebra $C$---a direct sum of matrix algebras over $C[0,1]$ or $\mathbb{C}$ as
$$\xi: A_{n}\xrightarrow{\xi_{1}}C\xrightarrow{\xi_{2}}QA_{m}Q.$$\\
The proof is similar to Proposition 2.8 and is omitted.

\noindent\textbf{2.10.} Let $\alpha:\mathbb{Z}\longrightarrow\mathbb{Z}/k_{1}\mathbb{Z}$ be the group homomorphism defined by $\alpha(\textbf{1})=[\textbf{1}]$, where the right hand side is the equivalent class $[\textbf{1}]$ of $\textbf{1}$ in $\mathbb{Z}/k_{1}\mathbb{Z}$. Then it is well known from homological algebra that for the group $\mathbb{Z}/k\mathbb{Z}$, $\alpha$ induces a surjective map
$$\alpha_{\ast}: Ext(\mathbb{Z}/k\mathbb{Z},\mathbb{Z})(=\mathbb{Z}/k\mathbb{Z})\longrightarrow Ext(\mathbb{Z}/k\mathbb{Z},\mathbb{Z}/k_{1}\mathbb{Z})(=\mathbb{Z}/(k,k_{1})\mathbb{Z}),$$
where $(k,k_{1})$ is the greatest common factor of $k$ and $k_{1}$.\\
Recall, as in [DN], for two connected finite simplicial complexes $X$ and $Y$, we use $kk(Y,X)$ to denote the group of equivalent classes of homomorphisms from  $C_{0}(X\backslash\{pt\})$ to $C_{0}(Y\backslash\{pt\})\otimes {\mathcal K}(H)$. Please see [DN] for details.

\noindent\textbf{Lemma 2.11.} (a) Any unital homomorphism
$$\phi: C(T_{\uppercase\expandafter{\romannumeral2}, k})\longrightarrow M_{\bullet}(C(T_{III,k_{1}})),$$
is homotopy equivalent to unital homomorphism $\psi$ factor as
$$C(T_{\uppercase\expandafter{\romannumeral2}, k})\xrightarrow{\psi_{1}}C(S^{1})\xrightarrow{\psi_{2}}M_{\bullet}(C(T_{III,k_{1}})).$$
(b) Any unital homomorphism $\phi: C(T_{\uppercase\expandafter{\romannumeral2}, k})\longrightarrow PM_{\bullet}(C(S^{2}))P$ is homotopy equivalent to unital homomorphism $\psi$ factor as
$$C(T_{\uppercase\expandafter{\romannumeral2}, k})\xrightarrow{\psi_{1}}\mathbb{C}\xrightarrow{\psi_{2}}PM_{\bullet}(C(S^{2}))P.$$
\begin{proof}
Part (b) is well known (see chapter 3 of [EG2]). To prove part (a), we note that
$$KK(C_{0}(S^{1}\backslash \{1\}),C_{0}(T_{III,k_{1}}\backslash \{x_{1}\})=kk(T_{III,k_{1}},S^{1})=\mathbb{Z}/k_{1}\mathbb{Z}=Hom(\mathbb{Z},\mathbb{Z}/k_{1}\mathbb{Z})),$$
where $x_{1}\in T_{III,k_{1}}$ is a base point. The map $\alpha: \mathbb{Z}\longrightarrow\mathbb{Z}/k_{1}\mathbb{Z}$ in 2.10 can be induced by a homomorphism: $\psi_{2}: C(S^{1})\longrightarrow M_{\bullet}(C(T_{\uppercase\expandafter{\romannumeral3}, k}))$.\\
Let
$$[\phi]\in kk(T_{III,k_{1}}\; T_{\uppercase\expandafter{\romannumeral2}, k})=Ext(K_{0}(C_{0}(T_{\uppercase\expandafter{\romannumeral2}, k}\backslash\{x_{0}\})), K_{1}(C_{0}(T_{\uppercase\expandafter{\romannumeral3}, k}))),$$
be the element induced by homomorphism $\phi$, where $\{x_{0}\}$ is the base point. By 2.10
$$[\phi]=\beta\times[\psi_{2}],\;for\;\beta\in kk(S^{1},T_{\uppercase\expandafter{\romannumeral2}, k})=Ext(K_{0}(C(T_{\uppercase\expandafter{\romannumeral2}, k}\backslash\{x_{0}\})),K_{1}(C(S^{1}))),$$
on the other hand $\beta$ can be realized by unital homomorphism
$$\psi_{1}: C(T_{\uppercase\expandafter{\romannumeral2}, k})\longrightarrow C(S^{1}).$$
(see section 3 of [EG2]).\\
\end{proof}
The following result is modification of Theorem 3.12 of [GJLP2].

\noindent\textbf{Theorem 2.12.} Let $B_{1}=\bigoplus_{i=1}^{s}B_{1}^{i}$, each $B^{i}$ is either matrix algebras over $\{pt\}, [0,1], S^{1}$ or $\{T_{\uppercase\expandafter{\romannumeral2}, k}\}^{\infty}_{k=2}$ or a dimension drop algebra.
Let $\varepsilon_{1}>0$ and let
$$\widetilde{G}_{1}(=\oplus\widetilde{G}^{i}_{1})\subset G_{1}(=\oplus G^{i}_{1})\subset B_{1}(=\oplus B^{i}_{1}),$$
be a finite set with $\omega(G^{i}_{1})<\varepsilon_{1}$ for $B^{i}_{1}=M_{\bullet}(C(T_{\uppercase\expandafter{\romannumeral2}, k}))$.

Let $A=M_{N}(C(X))$, where $X$ is one of form $\{pt\}, [0,1], S^{1}, \{T_{\uppercase\expandafter{\romannumeral2}, k}\}^{\infty}_{k=2}, \{T_{\uppercase\expandafter{\romannumeral3}, k}\}^{\infty}_{k=2}$ or $S^{2}$. Let $\alpha_{1}: B_{1}\longrightarrow A$ be a homomorphism. Let $F_{1}\subset A$ be a finite set, and let $\varepsilon(<\varepsilon_{1})$ and $\delta$ be any positive number. Then there exists a diagram.
$$
\xymatrix{
 A\ar[rdrd]^{\beta}\ar[rr]^{\phi} &  &  A^{'}  \\
     &  &           &  &         &  &          \\
     B_{1}\ar[uu]_{\alpha_{1}}\ar[rr]^{\psi}  &  &  B_{2}\ar[uu]_{\alpha_{2}}     \\
 }
$$
where $A^{'}=M_{K}(A)$, $B_{2}$ is described as below. $B_{2}$ is a direct sum of finite dimensional $C^{*}$-algebra and dimension drop algebra if $X=T_{\uppercase\expandafter{\romannumeral3}, k}$. $B_{2}$ is a finite dimensional algebra if $X=S^{2}$. And $B_{2}=M_{\bullet}(A)$ if $X=\{pt\}, [0,1], S^{1}, T_{\uppercase\expandafter{\romannumeral2}, k}$. Furthermore the diagram satisfies
the following conditions\\
(1) $\psi$ is a homomorphism, $\alpha_{2}$ is a unital injective homomorphism and $\phi$ is a unital simple embedding;\\
(2) $\beta\in Map(A,B_{2})_{1}$ is $F_{1}-\delta$ multiplicative;\\
(3) if $B^{i}_{1}$ is of form $M_{\bullet}(C(T_{\uppercase\expandafter{\romannumeral2}, k}))$, then
$$\parallel\psi(g)-\beta\circ\alpha_{1}(g)\parallel<10\varepsilon_{1},\; ~~\forall g\in\widetilde{G}^{i}_{1};$$
and if $B^{i}_{1}$ is not of form $M_{\bullet}(C(T_{\uppercase\expandafter{\romannumeral2}, k}))$, then
$$\parallel \psi(g)-\beta\circ\alpha_1(g)\parallel<\varepsilon,\; ~~\forall g\in G^{i}_{1};$$
(4) If $X=T_{\uppercase\expandafter{\romannumeral2}, k}$, then $\omega(\beta(F_{1})\cup\psi(G_{1}))<\varepsilon.$\\
(Note that, we only required that the weak variation of finite set in $M_{\bullet}(C(T_{\uppercase\expandafter{\romannumeral2}, k})$ to be small. In particular, we do not need to introduce the concept of weak variation for a finite subset of a dimension drop algebra.)
\begin{proof}.
For $X=T_{\uppercase\expandafter{\romannumeral2}, k},\{pt\}, [0,1]$ or $S^{1}$, one can choose $B_{2}=M_{K}(A)=A^{'}$ and homomorphism $\phi=\beta: A\longrightarrow B_{2}$ being simple embedding such that
$$\omega(\beta(F_{1})\cup\alpha_{1}(G_{1})))<\varepsilon.$$
This can be done by choosing $K$ large enough.
Choose $\psi=\beta\circ\alpha_{1}$, and $\alpha_{2}=id: B_{2}\longrightarrow A^{'}$.\\
For the case $X=T_{\uppercase\expandafter{\romannumeral3}, k}$, or $S^{2}$, the requirement (4) is an empty requirement.

We will deal with each block of $B_{1}$ separately. For the block $B^{i}_{1}$ other than $M_{\bullet}(C(T_{\uppercase\expandafter{\romannumeral2}, k})$, the construction can be done easily by using Lemma 2.7, since $B^{i}_{1}$ is stably generated, which implies that any sufficiently multiplicative map from $B^{i}_{1}$ is close to a homomorphism. So we assume that $B^{i}_{1}=M_{\bullet}(C(T_{\uppercase\expandafter{\romannumeral2}, k})$. Recall we already assume $A$ is of form $M_{\bullet}(C(T_{\uppercase\expandafter{\romannumeral3}, k}))$ or $M_{\bullet}(C(S^{2})$. By Lemma 2.11, the homomorphism $\alpha_{1}: B^{i}_{1}\longrightarrow A$ is a homotopy to $\alpha^{'}: B^{i}_{1}\longrightarrow A$ with $\alpha^{'}(\textbf{1}_{B^{i}_{1}})=\alpha_{1}(\textbf{1}_{B^{i}_{1}})$ and $\alpha^{'}$ factor as
$$B_{1}^i\xrightarrow{\xi_{1}}C\xrightarrow{\xi_{2}}A,$$
where $C$ is a finite dimensional $C^{*}$-algebra for the case $X=S^{2}$ or $C=M_{\bullet}(C(S^1)$  for the case $X=T_{\uppercase\expandafter{\romannumeral3}, k}$ (note that $B^{i}_{1}=M_{\bullet}(C(T_{\uppercase\expandafter{\romannumeral2}, k})$). Since $C$ is stably generated, there is a finite set $E_{1}\subset A$ and $\delta_{1}>0$ such that if a complete positive map $\beta: A\longrightarrow D$ (for any $C^{*}$-algebra $D$) is $E_{1}-\delta_{1}$ multiplicative, then the map $\beta\circ\xi_{2}: C\longrightarrow D$ can be perturbed to a homomorphism
$$\widetilde{\xi}: C\longrightarrow D$$
such that
$$\|\widetilde{\xi}(g)-\beta\circ\xi_{2}(g)\|<\varepsilon_{1},\;\forall\;g\in\xi_{1}(\widetilde{G}^{i}_{1}).$$
Apply Theorem 1.6.9 of [G5] to two homotopic homomorphism
$$\alpha_{1},\alpha^{'}: B^{i}_{1}\longrightarrow A,$$
and $G^{i}_{1}\subset B^{i}_{1}$ which is approximately constant to within $\varepsilon_{1}$, to obtain a finite set $E_{2}\subset A$, $\delta_{2}>0$ and positive integer $L^{'}>0$ (in places of $G,\delta$ and $L$ in Theorem 1.6.9 of [G5]). Apply Lemma 2.7 to the set $\widetilde{E}=E_{1}\cup E_{2}\cup F_{1}$ and $\widetilde{\delta}=\frac{1}{3}\min(\varepsilon,\delta,\delta_{1},\delta_{2})$ to obtain the diagram
$$
\xymatrix{
A\ar[rr]^{\phi^{'}}\ar[rdrd]^{\beta^{'}}  &  &  M_{L_{1}}(A)  &  &  \\
     &  &  &  &  &  &  \\
                                                                     &  &   B^{'}\ar[uu]_{\iota} & &   \\
 }
$$
with $\beta^{'}$ being $\widetilde{E}-\widetilde{\delta}$ multiplicative and
$$\parallel\iota\circ\beta^{'}(f)-\phi^{'}(f)\parallel<\widetilde{\delta},\;\forall\;f\in\widetilde{E}.$$
Let $L=L^{'}\cdot rank(\textbf{1}_{A})$, and let $\beta_{1}: A\longrightarrow M_{L}(B^{'})$ be any unital homomorphism defined by point evaluation. Then by Theorem 1.6.9 of [G5], there is a unitary $u\in M_{L+1}(B)$ such that
$$\|u((\beta^{'}\oplus\beta_{1})\circ\alpha^{'}(f))u^{\ast}-(\beta^{'}\oplus\beta_{1})\circ\alpha_{1}(f)\|<8\varepsilon_{1},\;\forall\;f\in\widetilde{G}^{i}_{1}.$$ By the choice of $E_{1}$, there is a homomorphism
$$\widetilde{\xi}: C\longrightarrow M_{L+1}(B^{'}),$$
such that
$$\|\widetilde{\xi}(f)-u((\beta^{'}\oplus\beta)\circ\xi_{2}(f))u^{\ast}\|<\varepsilon_{1},\;\forall\;f\in\xi_{1}(\widetilde{G}^{i}_{1}).$$ Define $B_{2}=M_{L+1}(B^{'}), K=L_{1}(L+1), A^{'}=M_{K}(A)=M_{L+1}(M_{L_{1}}(A))$,
\begin{center}
$\psi: B^{i}_{1}\longrightarrow B_{2}$~~by~~$\psi=\widetilde{\xi}\circ\xi_{1}: B^{i}_{1}\xrightarrow{\xi_{1}}C\xrightarrow{\widetilde{\xi}}B_{2}$,
\end{center}
\begin{center}
$\beta: A\longrightarrow M_{L+1}(B^{'})$~~by~~$\beta=\beta^{'}\oplus\beta_{1}$,
\end{center}
and
\begin{center}
$\phi: A\longrightarrow M_{L+1}(M_{L_{1}}(A))$~~by~~$\phi=\phi^{'}\oplus((\iota\otimes id_{L})\circ\beta_{1})$,
\end{center}
(note that $\beta_{1}$ is a homomorphism) to finish the proof.\\
\end{proof}

\noindent\textbf{2.13.} Recall for $A=\bigoplus^{t}_{i=1}M_{k_{i}}(C(X_{i}))$, where $X_{i}$ are  path connected simplicial complexs, we use the notation $r(A)$ to denote $\bigoplus^{t}_{i=1}M_{k_{i}}(\mathbb{C})$ which could be considered to be the subalgebra consisting of all t-tuples of constant function from $X_{i}$ to $M_{k_{i}}(\mathbb{C})$ ($i=1,2,\cdots,t$). Fixed a base point $x^{0}_{i}\in X_{i}$ for each $X_{i}$, one defines a map $r: A\longrightarrow r(A)$ by
$$r(f_{1},f_{2},\cdots,f_{t})=(f_{1}(x^{0}_{1}),f_{2}(x^{0}_{2}),\cdots,f_{t}(x^{0}_{t}))\in r(A).$$
We have the following Corollary

\noindent\textbf{Corollary 2.14.} Let $B_{1}=\oplus B^{j}_{1}$, where $B^{j}_{1}$ are either of form $M_{k(j)}(C(X_{j}))$, with $X_{j}$ being one of
$\{pt\}$, [0,1], $S^{1}$, $\{T_{\uppercase\expandafter{\romannumeral2}, k}\}^{\infty}_{k=2}$ or one of form $M_{k(j)}(I_{l(j)})$. Let $\alpha_{1}; B_{1}\longrightarrow A$ be a homomorphism, where $A$ is a direct sum of matrix algebras over $\{pt\}, [0,1], S^{1}, \{T_{\uppercase\expandafter{\romannumeral2}, k}\}^{\infty}_{k=2}, \{T_{\uppercase\expandafter{\romannumeral3}, k}\}^{\infty}_{k=2}$ and $S^{2}$.
Let $\varepsilon_{1}>0$ and let $\widetilde{E}(=\oplus\widetilde{E}^{i})\subset E(=\oplus E^{i})\subset B_{1}(=\oplus B^{i}_{1})$ be two finite subset with the condition
\begin{center}
$\omega(\widetilde{E}^{i})<\varepsilon_{1}$, if $B^{i}_{1}=M_{\bullet}(C(Y_{i}))$ with $Y_{i}\in\{T_{\uppercase\expandafter{\romannumeral2}, k}\}^{\infty}_{k=2}$.
\end{center}
Let $F\subset A$ be any finite set, $\varepsilon_{2}>0$, $\delta>0$. Then there exists a diagram
$$
\xymatrix{
 A\ar[rdrd]^{\beta\oplus r}\ar[rr]^{\phi\oplus r} &  &  A^{'}\oplus r(A)  \\
     &  &           &  &         &  &          \\
     B_{1}\ar[uu]_{\alpha_{1}}\ar[rr]^{\psi\oplus(r\circ\alpha_{1})}  &  &  B_{2}\oplus r(A)\ar[uu]_{\alpha_{2}\oplus id}     \\
 }
$$
where $A^{'}=M_{L}(A)$, and $B_{2}$ is a direct sum of matrix algebras over space: $\{pt\}$, [0,1], $S^{1}$, $\{T_{\uppercase\expandafter{\romannumeral2}, k}\}^{\infty}_{k=2}$ and dimension drop algebras, with the following properties.  \\
(1) $\psi$ is a homomorphism, $\alpha_{2}$ is a injective homomorphism and $\phi$ is a unital simple embedding.\\
(2) $\beta\in Map(A,B_{2})_{1}$ is $F_{1}-\delta$ multiplicative.\\
(3) for $g\in\widetilde{E}^{i}$ with $B^{i}_{1}=M_{\bullet}(C(X_{i}))$, $X_{i}\in\{T_{\uppercase\expandafter{\romannumeral2}, k}\}^{\infty}_{k=2}$ we have
$$\|(\beta\oplus r)(g)-(\psi\oplus(r\circ\alpha_{1}))(g)\|\leqslant10\varepsilon_{1};$$
and for $g\in E^{i}(\supset\widetilde{E}^{i})$ where $B^{i}_{1}$ is not of form $M_{\bullet}(C(T_{\uppercase\expandafter{\romannumeral2}, k}))$, we have
$$\|(\beta\oplus r)(g)-(\psi\oplus(r\circ\alpha_{1}))(g)\|<\varepsilon_{1};$$
and for $f\in F$
$$\|(\alpha_{2}\oplus id)\circ(\beta\oplus r)(f)-(\phi\oplus r)(f)\|<\varepsilon_{1}.$$
(4) For $B^{i}_{2}$ of form $M_{\bullet}(C(T_{\uppercase\expandafter{\romannumeral2}, k}))$,
$$\omega(\pi_{i}(\beta(F)\cup\psi(E)))<\varepsilon_{2},$$
where $\pi_{i}$ is the canonical projection from $B_{2}$ to $B^{i}_{2}$.\\

{\bf Remark}: In the application of the above Corollary, we will denote the map $\beta\oplus r$ by $\beta$ and $\psi\oplus(r\circ\alpha_{1})$ by $\psi$.

\vspace{3mm}

\noindent\textbf{\S3. The proof of main theorem}

\vspace{-2.4mm}

In this section, we prove the following main theorem.

\noindent\textbf{Theorem 3.1.} Suppose $\lim(A_{n}=\bigoplus_{i=1}^{t_{n}}M_{[n,i]}(C(X_{n,i})), \phi_{n,m})$ is an $AH$ inductive limit with $X_{n,i}$ being among the spaces $\{pt\}, [0,1], S^{1}, \{T_{\uppercase\expandafter{\romannumeral2}, k}\}^{\infty}_{k=2}$, $\{T_{\uppercase\expandafter{\romannumeral3}, k}\}^{\infty}_{k=2}$, such that the limit algebra $A$ has ideal property. Then there is another inductive system ($B_{n}=\oplus B^{i}_{n}, \psi_{n,m}$) with same limit algebra, where $B^{i}_{n}$ are either $M_{[n,i]^{'}}(C(Y_{n,i}))$ with $Y_{n,i}$ being one of $\{pt\}, [0,1], S^{1}, \{T_{\uppercase\expandafter{\romannumeral2}, k}\}^{\infty}_{k=2}$ (but without $T_{\uppercase\expandafter{\romannumeral3}, k}$ and $S^{2}$), or dimension drop algebra $M_{[n,i]^{'}}(I_{k(n,i)})$.
 \begin{proof}
 Let $\varepsilon_{1}>\varepsilon_{2}>\varepsilon_{3}>\cdots$ be a sequence of  positive number with $\sum\varepsilon_{n}<+\infty$. We need to construct the intertwining diagram

%
 $$
 \xymatrix@R=0.5ex{
 F_{1} & & F_{2} & & & & F_{n} & & F_{n+1}\\
 \bigcap & & \bigcap & & & & \bigcap & & \bigcap\\
 A_{s(1)} \ar[rdrd]^{\beta_{1}} \ar[rr]^{\phi_{s(1),s(2)}} & & A_{s(2)} \ar[rdrd]^{\beta_{2}} \ar[rr]^{\phi_{s(2),s(3)}} & &\cdots \ar[rr] & & A_{s(n)} \ar[rdrd]^{\beta_{n}} \ar[rr]^{\phi_{s(n),s(n+1)}}  & & A_{s(n+1)} \ar[rr] \ar[rdrd] & & \cdots\\
   & & & & & &\\
 B_{1} \ar[rr]^{\psi_{1,2}} \ar[uu]^{\alpha_{1}} & & B_{2} \ar[uu]^{\alpha_{2}} \ar[rr]^{\psi_{2,3}} & &\cdots \ar[rr] & & B_{n}\ar[rr]^{\psi_{n,n+1}} \ar[uu]^{\alpha_{n}} &  & B_{n+1} \ar[uu]^{\alpha_{n+1}} \ar[rr] & & \cdots\\
 \bigcup & & \bigcup & & & & \bigcup & & \bigcup\\
 E_{1} & & E_{2} & & & & E_{n} & & E_{n+1}\\
 \bigcup & & \bigcup & & & & \bigcup & & \bigcup\\
 \widetilde{E_{1}} & & \widetilde{E_{2}} & & & & \widetilde{E_{n}} & & \widetilde{E_{n+1}}
 }
 $$
satisfying the following conditions\\
(0.1) $(A_{s(n)},\phi_{s(n),s(m)})$ is a sub-inductive system of $(A_{n}, \phi_{n,m}), (B_{n}, \psi_{n,m})$ is an inductive system of direct sum of matrix algebras
over the spaces $\{pt\}, [0,1], S^{1}, {T_{\uppercase\expandafter{\romannumeral2},k}}$ and dimension drop algebra $M_{\bullet}(I_{k(n,i)})$.\\
(0.2) Choose $\{a_{i,j}\}^{\infty}_{j=1}\subset A_{s(i)}$ and $\{b_{i,j}\}^{\infty}_{j=1}\subset B_{i}$ to be countable dense subsets of unit balls of $A_{s(i)}$ and $B_{i}$, respectively. $F_{n}$ are subsets of unit balls of $A_{s(n)}$, and $\widetilde{E_{n}}\subset E_{n}$ are both subsets of unit balls of $B_{n}$ satisfying
$$\phi_{s(n),s(n+1)}(F_{n})\cup\alpha_{n+1}(E_{n+1})\cup\bigcup^{n+1}_{i=1}\phi_{s(i),s(n+1)}(\{a_{i1},a_{i2},\cdot\cdot\cdot, a_{in+1}\})\subset F_{n+1},$$
$$\psi_{n,n+1}(E_{n})\cup\beta_{n}(F_{n})\subset \widetilde{E}_{n+1}\subset E_{n+1},$$
and
$$\bigcup^{n+1}_{i=1}\psi_{i,n+1}(\{b_{i1},b_{i2},\cdot\cdot\cdot, b_{in+1}\})\subset E_{n+1}.$$
(Here $\phi_{n,n}:A_{n}\longrightarrow A_{n}$, and $\psi_{n,n}:B_{n}\longrightarrow B_{n}$ are understand as identity maps.)\\
(0.3) $\beta_{n}$ are $F_{n}-2\varepsilon_{n}$ multiplicative and $\alpha_{n}$ are homomorphism.\\
(0.4) for all $g\in\widetilde{E}_{n}$, $$\|\psi_{n,n+1}(g)-\beta_{n}\circ\alpha_{n}(g)\| < 14\varepsilon_{n},$$
and for all $f\in F_{n}$,
 $$\|\phi_{s(n),s(n+1)}(f)-\alpha_{n+1}\circ\beta_{n}(f)\| < 14\varepsilon_{n}.$$
(0.5) For any block $B^{i}_{n}$ with spectrum $T_{\uppercase\expandafter{\romannumeral2},k},$ we have $\omega(\widetilde{E}_{n}^{i}) < \varepsilon_{n}$, where $\widetilde{E}_{n}^{i}=\pi_{i}(\widetilde{E}_{n})$ for $\pi_{i}:B_{n}\longrightarrow B^{i}_{n}$ the canonical projections.

The diagram will be constructed inductively. First, let $B_{1}=\{0\}, A_{s(1)}=A_{1}, \alpha_{1}=0$. Let $b_{1j}=0\in B_{1}$ for $j=1,2,...$ and let $\{a_{1j}\}^{\infty}_{j=1}$ be a countable dense subset of the unit ball of $A_{s(1)}$. And let $\widetilde{E}_{1}=E_{1}=\{b_{11}\}=B_{1}$
and $F_{1}=\bigoplus^{t_{1}}_{i=1}F^{i}_{1}$, where $F^{i}_{1}=\pi_{i}(\{a_{11}\})\subset A^{i}_{1}$.

As inductive assumption, assume that we already have the diagram

%

 $$
 \xymatrix@R=0.4ex{
F_{1} & & F_{2} & & & & F_{n}\\
\bigcap & & \bigcap & & & & \bigcap\\
 A_{s(1)} \ar[rdrd]^{\beta_{1}} \ar[rr]^{\phi_{s(1),s(2)}} & & A_{s(2)} \ar[rdrd]^{\beta_{2}} \ar[rr]^{\phi_{s(2),s(3)}} & &\cdots \ar[rr] \ar[rdrd]^{\beta_{n-1}} & & A_{s(n)}\\
   & & & & & &\\
 B_{1} \ar[rr]^{\psi_{1,2}} \ar[uu]^{\alpha_{1}} & & B_{2} \ar[uu]^{\alpha_{2}} \ar[rr]^{\psi_{2,3}} & &\cdots \ar[rr] & & B_{n}\ar[uu]^{\alpha_{n}}\\
 \bigcup & & \bigcup & & & & \bigcup\\
 E_{1} & & E_{2} & & & & E_{n}\\
 \bigcup & & \bigcup & & & & \bigcup\\
 \widetilde{E_{1}} & & \widetilde{E_{2}} & & & & \widetilde{E_{n}}
 }
 $$

and for each $i=1,2,\cdot\cdot\cdot,n$, we have dense subsets
$\{a_{ij}\}^{\infty}_{j=1}\subset$ the unit ball of $A_{s(i)}$
and
$\{b_{ij}\}^{\infty}_{j=1}\subset$ the unit ball of $B_{i}$
satisfying  the conditions (0.1)-(0.5) above. We have to construct the next piece of the diagram
$$
\xymatrix@!C{
F_{n}\subset A_{s(n)}\ar[rr]^{\phi_{s(n),s(n+1)}}\ar[rdrd]^{\beta_{n}}  &  &  A_{s(n+1)}\supset F_{n+1}  &  &    \\
     &  &  &  &  &  &  \\
\widetilde{E}_{n}\subset E_{n}
\subset B_{n}\ar@<-8mm>[uu]^{\alpha_{n}}\ar[rr]_{\psi_{n,n+1}} &  &
   \;\;B_{n+1}\ar@<10mm>[uu]_{\alpha_{n+1}}\supset
E_{n+1}\supset\widetilde{E}_{n+1}  &  &
 }
$$
to satisfy the condition (0.1)-(0.5).

Among the conditions for induction assumption, we will only use the conditions that $\alpha_{n}$ is a homomorphism and (0.5) above.\\
\textbf{Step 1}. We enlarge $\widetilde{E}_{n}$ to $\bigoplus_{i}\pi_{i}(\widetilde{E}_{n}^{i})$ and $E_{n}$ to $\bigoplus_{i}\pi_{i}(E_{n})$. Then $\widetilde{E}_{n}(=\oplus\widetilde{E}_{n}^{i})\subset E_{n}(=\oplus E_{n})$ and for each $B^{i}_{n}$ with spectrum $T_{\uppercase\expandafter{\romannumeral2},k}$, we have $\omega(E_{n}^{i})<\varepsilon_{n}$ from induction assumption (0.5). By proposition 2.9 applied to $\alpha_{n}:B_{n}\rightarrow A_{s(n)},\widetilde{E}_{n}\subset E_{n}\subset B_{n},F_{n}\subset A_{s(n)}$ and $\varepsilon_{n}>0$, there are $A_{m_{1}}(m_{1}>s(n))$, two othogonal projections $P_{0},P_{1}\in A_{m_{1}}$ with $\phi_{s(n),m_{1}}(\textbf{1}_{A_{s(n)}})=P_{0}+P_{1}$ and $P_{0}$ trivial, a $C^{*}$-algebra $C$---a direct sum of matrix algebras over $C[0,1]$ or $\mathbb{C}$, a unital map $\theta\in Map(A_{s(n)},P_{0}A_{m_{1}}P_{0})_{1}$, a unital homomorphism $\xi_{1}\in Hom(A_{s(n)},C)_{1}$, a unital homomorphism $\xi_{2}\in Hom(C,P_{1}A_{m_{1}}P_{1})_{1}$ such that\\
(1.1) $\|\phi_{s(n),m_{1}}(f)-\theta(f)\oplus(\xi_{2}\circ\xi_{1})(f)\|<\varepsilon_{n}$ for all $f\in F_{n}$.\\
(1.2) $\theta$ is $F_{n}-\varepsilon$ multiplicative and $F:=\theta(F_{n})$ satisfy $\omega(F)<\varepsilon_{n}$.\\
(1.3) $\|\alpha(g)-\theta\circ\alpha_{n}(g)\|<3\varepsilon_{n}$ for all $g\in\widetilde{E}_{n}$.

Let all the blocks of C be parts of $C^{\ast}$-algebra $B_{n+1}$. That is
\begin{center}
$B_{n+1}=C\oplus$ (some other blocks).
\end{center}
The map $\beta_{n}:A_{s(n)}\rightarrow B_{n+1}$, and the homomorphism $\psi_{n,n+1}:B_{n}\rightarrow B_{n+1}$ are defined by $\beta_{n}=\xi_{1}:A_{s(n)}\rightarrow C(\subset B_{n+1})$
and $\psi_{n,n+1}=\xi_{1}\circ\alpha_{n}:B_{n}\rightarrow C(\subset B_{n+1})$
for the blocks of $C(\subset B_{n+1})$. For this part, $\beta_{n}$ is also a homomorphism.\\
\textbf{Step 2}. Let $A=P_{0}A_{m_{1}}P_{0},F=\theta(F_{n})$. Since $P_{0}$ is a trivial projection, $$A\cong\oplus M_{l_{i}}(C(X_{m_{1},i})).$$ Let $r(A):=\oplus M_{l_{i}}(\mathbb{C})$ and $r:A\rightarrow r(A)$ be as in 2.13. Applying corollary 2.14 and its remark to $\alpha:B_{n}\rightarrow A,\widetilde{E}_{n}\subset E_{n}\subset B_{n}$ and $F\subset A$, we obtain the following diagram
$$
\xymatrix{
A\ar[rr]^{\phi\oplus r}\ar[rdrd]^{\beta}      & &M_{L}(A)\oplus r(A)  &  &    \\
     &  &  &  &  &  &  \\
B_{n}\ar[rr]^{\psi}\ar[uu]^{\alpha} &  &  B\ar[uu]_{\alpha^{'}} &  &  \\
 }
$$
such that\\
(2.1) $B$ is a direct sum of matrix algebras over $\{pt\},[0,1],S^{1},T_{\uppercase\expandafter{\romannumeral2},k}$ and dimension drop algebra.\\
(2.2) $\alpha^{'}$ is an injective homomorphism and $\beta$ is $F-\varepsilon_{n}$ multiplicative.\\
(2.3) $\phi:A\rightarrow M_{L}(A)$ is a unital simple embedding and $r:A\rightarrow r(A)$ is as in 2.13.\\
(2.4) $\|\beta\circ\alpha(g)-\psi(g)\|<10\varepsilon_{n}$ for all $g\in \widetilde{E}_{n}$ and $\|(\phi\oplus r)(f)-\alpha^{'}\circ\beta(f)\|<\varepsilon_{n}$ for all $f\in F(:=\theta(F_{n}))$.\\
(2.5) $\omega(\pi_{i}(\psi(E_{n}))\cup\beta(F))<\varepsilon_{n+1}$ (note that $\beta(F)=\beta\circ\theta(F_{n})$), for $B_n^{i}$ being of form $M_{\bullet}(C(X))$ with $X\in\{T_{\uppercase\expandafter{\romannumeral2}, k}\}^{\infty}_{k=2}$.

Let all the blocks $B$ be also part of $B_{n+1}$, that is
\begin{center}
$B_{n+1}=C\oplus B\oplus$ (some other blocks).
\end{center}
The maps $\beta_{n}:A_{s(n)}\longrightarrow B_{n+1},\psi_{n,n+1}:B_{n}\longrightarrow B_{n+1}$
are defined by
$$\beta_{n}:=\beta\circ\theta:A_{s(n)}\xrightarrow{\theta}A\xrightarrow{\beta}B(\subset B_{n+1}),$$
and
$$\psi_{n,n+1}:=\psi:B_{n}\rightarrow B(\subset B_{n+1}),$$ for the blocks of $B(\subset B_{n+1})$. This part of $\beta_{n}$ is $F_{n}-2\varepsilon_{n}$ multiplicative, since $\theta$ is $F_{n}-\varepsilon_{n}$ multiplicative, $\beta$ is $F-\varepsilon_{n}$ multiplicative and $F=\theta(F_{n})$.\\
\textbf{Step 3}. By Lemma 3.15 of [GJLP2] applied to $\phi\oplus r:A\rightarrow M_{L}(A)\oplus r(A)$, there is an $A_{m_{2}}$ and there is a unital homomorphism
$$\lambda:M_{L}(A)\oplus r(A)\rightarrow RA_{m_{2}}R,$$ where $R=\phi_{m_{1},m_{2}}(P_{0})$ (write $R$ as $\bigoplus_{j}R^{j}\in\bigoplus_{j}A^{j}_{m}$) such that the diagram
$$
\xymatrix{
A(=P_{0}A_{m}P_{0})\ar[rr]^{\phi_{m_{1},m_{2}}}\ar[rdrd]^{\phi\oplus r}  &  &  RA_{m_{2}}R  &  &  \\
     &  &  &  &  &  &  \\
       &  &   M_{L}(A)\oplus r(A)\ar[uu]_{\lambda} & &
 }
$$
satisfies the following conditions:\\
(3.1) $\lambda\circ(\phi\oplus r)$ is homotopy equivalent to
$$\phi^{'}:=\phi_{m_{1},m_{2}}|_{A}.$$
\textbf{Step 4}. Applying Theorem 1.6.9 of [G5] to finite set $F\subset A$ (with $\omega(F)<\varepsilon_{n}$)
and to two homotopic homomorphisms $\phi^{'}$ and $\lambda\circ(\phi\oplus r):A\rightarrow RA_{m_{2}}R$ (with $RA_{m_{2}}R$ in place of $C$ in Theorem 1.6.9 of [G5]), we obtain a finite set $F^{'}\subset RA_{m_{2}}R,\; \delta>0$ and $L>0$ as in the Theorem.

Let $G=\oplus\pi_{i}(\psi(E_{n})\cup\beta(F))=\oplus G^{i}$. Then by (2.5), we have $\omega(G^{i})<\varepsilon_{n+1}$, if $B^{i}$ is of form $M_{\bullet}(C(T_{\uppercase\expandafter{\romannumeral2}, k}))$. By Proposition 2.8 applied to $RA_{m_{2}}R$ and
$$\lambda\circ\alpha^{'}:B\rightarrow RA_{m_{2}}R,$$ finite set $G\subset B$, $F^{'}\cup(\phi_{m_{1}m_{2}}\mid_{A}(F))\in RA_{m_{2}}R$, $min(\varepsilon_{n},\delta)>0$ (in place of $\varepsilon$) and $L>0$, there are $A_{s(n+1)}$, mutually orthogonal projections $Q_{0},Q_{1},Q_{2}\in A_{s(n+1)}$ with $\phi_{m_{2},s(n+1)}(R)=Q_{0}\oplus Q_{1}\oplus Q_{2}$, a $C^{*}$-algebra $D$---a direct sum of matrix algebras over C[0,1] or $\mathbb{C}$---,
a unital map $\theta_{0}\in$ Map($RA_{m_{2}}R,Q_{0}A_{s(n+1)}Q_{0}$), and four unital homomorphisms $$\theta_{1}\in Hom(RA_{m_{2}}R,Q_{1}A_{s(n+1)}Q_{1})_{1},\xi_{3}\in Hom(RA_{m_{2}}R,D)_{1},\xi_{4}\in Hom(D,Q_{2}A_{s(n+1)}Q_{2})_{1}$$ and $\alpha^{''}\in Hom(B,(Q_{0}+Q_{1})A_{s(n+1)}(Q_{0}+Q_{1}))_{1}$ such that the following is true.\\
(4.1) $\|\phi_{m_{2},s(n+1)}(f)-((\theta_{0}+\theta_{1})\oplus\xi_{4}\circ\xi_{3})(f)\|<\varepsilon_{n}$, for all $f\in\phi_{m_{1},m_{2}}|_{A}(F)\subset RA_{m_{2}}R$.\\
(4.2) $\|\alpha^{''}(g)-(\theta_{0}+\theta_{1})\circ\lambda\circ\alpha^{'}(g)\|<3\varepsilon_{n+1}<3\varepsilon_{n},\;\forall\;g\in G$.\\
(4.3) $\theta_{0}$ is $F^{'}-min(\varepsilon_{n},\delta)$ multiplicative and $\theta_{1}$ satisfies that
$$\theta^{i,j}_{1}([q])>L\cdot[\theta^{i,j}_{0}(R^{i})],$$
for any non zero projection $q\in R^{i}A_{m_{1}}R^{i}$.\\
By Theorem 1.6.9 of [G5], there is a unitary $u\in(Q_{0}\oplus Q_{1})A_{s(n+1)}(Q_{0}+Q_{1})$ such that $$\|(\theta_{0}+\theta_{1})\circ\phi^{'}(f)-Adu\circ(\theta_{0}+\theta_{1})\circ\lambda\circ(\phi\oplus r)(f)\|<8\varepsilon_{n},$$
for all $f\in F$.\\
Combining with second inequality of (2.4), we have\\
(4.4)$\|(\theta_{0}+\theta_{1})\circ\phi^{'}(f)-Adu\circ(\theta_{0}+\theta_{1})\circ\lambda\circ\alpha^{'}\circ\beta(f)\|<9\varepsilon_{n}$
for all $f\in F$.\\
\textbf{Step 5}. Finally let all blocks of $D$ be the rest of $B_{n+1}$. Namely, let $$B_{n+1}=C\oplus B\oplus D,$$ where $C$ is from Step 1, $B$ is from Step 2 and $D$ is from Step 4.

We already have the definition of $\beta_{n}:A_{s(n)}\rightarrow B_{n+1}$ and $\psi_{n,n+1}: B_{n}\rightarrow B_{n+1}$
for those blocks of $C\oplus B\subset B_{n+1}$
(from Step 1 and Step 2). The definition of $\beta_{n}$ and $\psi_{n,n+1}$ for blocks of $D$ and the homomorphism $\alpha_{n+1}:C\oplus B\oplus D\rightarrow A_{s(n+1)}$
will be given below.

The part of $\beta_{n}:A_{s(n)}\rightarrow D(\subset B_{n+1})$
is defined by
$$\beta_{n}=\xi_{3}\circ\phi^{'}\circ\theta:A_{s(n)}\xlongrightarrow{\theta}A\xlongrightarrow{\phi}RA_{m_{2}}R\xlongrightarrow{\xi_{3}}D.$$
(Recall that $A=P_{0}A_{m_{2}}P_{0}$ and $\phi^{'}=\phi_{m_{1},m_{2}}|_{A}.)$ Since $\theta$ is $F_{n}-\varepsilon_{n}$ multiplicative, and $\phi^{'}$ and $\xi_{3}$ are homomorphism, we know this part of $\beta_{n}$ is $F_{n}-\varepsilon_{n}$ multiplicative.

The part of $\psi_{n,n+1}:B_{n}\rightarrow D(\subset B_{n+1})$ is defined by
$$\psi_{n,n+1}=\xi_{3}\circ\phi^{'}\circ\alpha:B_{n}\xlongrightarrow{\alpha}A\xlongrightarrow{\phi^{'}}RA_{m}R\xlongrightarrow{\xi_{3}}D,$$
which is a homomorphism.

The homomorphism $\alpha_{n+1}:C\oplus B\oplus D\rightarrow A_{s(n+1)}$
is defined as following.

Let $\phi^{''}=\phi_{m_{1},s(n+1)}|_{P_{1}A_{m_{1}}P_{1}}:P_{1}A_{m_{1}}P_{1}\longrightarrow
\phi_{m_{1},s(n+1)}(P_{1})A_{s(n+1)}\phi_{m_{1},s(n+1)}(P_{1})$, where $P_{1}$ is from Step 1.
Define $$\alpha_{n+1}|_{C}=\phi^{''}\circ\xi_{2}:C\xlongrightarrow{\xi_{2}}P_{1}A_{m_{1}}P_{1}\xlongrightarrow{\phi^{''}}\phi_{m_{1},s(n+1)}(P_{1})A_{s(n+1)}\phi_{m_{1},s(n+1)}(P_{1}),$$
where $\xi_{2}$ is from Step 1.
\begin{center}
$\alpha_{n+1}|_{B}=Adu\circ\alpha^{''}:B\xlongrightarrow{\alpha^{''}}(Q_{0}\oplus Q_{1})A_{s(n+1)}(Q_{0}+Q_{1})\xlongrightarrow{Ad u}(Q_{0}\oplus Q_{1})A_{s(n+1)}(Q_{0}+Q_{1})$
\end{center}
where $\alpha^{''}$ is from Step 4, and define
$$\alpha_{n+1}|_{D}=\xi_{4}:D\rightarrow Q_{2}A_{s(n+1)}Q_{2}.$$
Finally choose $\{a_{n+1,j}\}^{\infty}_{j=1}\subset A_{s(n+1)}$ and $\{b_{n+1,j}\}^{\infty}_{j=1}\subset B_{n+1}$
to be countable dense subsets of the unit balls of $A_{s(n+1)}$ and $B_{n+1}$, respectively. And choose\\
$$F^{'}_{n+1}=\phi_{s(n),s(n+1)}(F_{n})\cup\alpha_{n+1}(E_{n+1})\cup\bigcup^{n+1}_{i=1}\phi_{s(i),s(n+1)}(\{a_{i1},a_{i2},\cdot\cdot\cdot, a_{in+1}\}),$$
$$E^{'}_{n+1}=\psi_{n,n+1}(E_{n})\cup\beta_{n}(F_{n})\cup\bigcup^{n+1}_{i=1}\psi_{i,n+1}(\{b_{i1},b_{i2},\cdot\cdot\cdot, b_{in+1}\}),$$
$$\widetilde{E}_{n+1}^{'}=\psi_{n,n+1}(E_{n})\cup\beta_{n}(F_{n})\subset E_{n+1}^{'}.$$
Define $F^{i}_{n+1}=\pi_{i}(F_{n+1}^{'})$ and $F_{n+1}=\bigoplus_{i}F^{i}_{n+1}$, $E^{i}_{n+1}=\pi_{i}(E^{'}_{n+1})$
and $E_{n+1}=\oplus_{i}E^{i}_{n+1}$.
For those block $B^{i}_{n+1}$ inside the algebra $B$ define $\widetilde{E}_{n+1}^{i}=\pi_{i}(\widetilde{E}_{n+1})$.
For those blocks inside $C$ and $D$, define $\widetilde{E}_{n+1}^{i}=E^{i}_{n+1}$.
And finally let $E_{n+1}=\bigoplus_{i}\widetilde{E}_{n+1}^{i}$.
Note all the blocks with spectrum $T_{\uppercase\expandafter{\romannumeral2},k}$ are in $B$. And (2.5) tells us that for those blocks $\omega(\widetilde{E}_{n+1}^{i})<\varepsilon_{n+1}$. Thus we obtain the following diagram
$$
\xymatrix@!C{
F_{n}\subset A_{s(n)}\ar[rr]^{\phi_{s(n),s(n+1)}}\ar[drdr]^{\beta_{n}}  &  &  A_{s(n+1)}\supset F_{n+1}  &  &    \\
     &  &  &  &  &  &  \\
\widetilde{E}_{n}\subset E_{n}
\subset B_{n}\ar@<-8mm>[uu]^{\alpha_{n}}\ar[rr]_{\psi_{n,n+1}} &  &
   \;\;B_{n}\ar@<12mm>[uu]_{\alpha_{n+1}}\supset
E_{n+1}\supset\widetilde{E}_{n+1}  &  &  \\
 }
$$
\textbf{Step 6}.Now we need to verify all the condition (0.1)-(0.5) for the above diagram.

From the end of Step 5, we know (0.5) holds, (0.1)-(0.2) hold from the construction (see the construction of $B, C, D$ in Step 1, 2 and 4, and $\widetilde{E}_{n+1}\subset E_{n+1}, F_{n+1}$ is the end of Step 5).
(0.3) follows from the end of Step 1, the end of Step 2 and the part of definition of $\beta_{n}$ for $D$ from Step 5.

So we only need to verify (0.4).

Combining (1.1) with (4.1), we have
$$\|\phi_{s(n),s(n+1)}(f)-[(\phi^{''}\circ\xi_{2}\circ\xi_{1})\oplus(\theta_{0}+\theta_{1})\circ\phi^{'}\circ\theta\oplus(\xi_{4}\circ\xi_{3}\circ\phi^{'}\circ\theta)(f)](f)\|
<\varepsilon_{n}+\varepsilon_{n}=2\varepsilon_{n}$$
for all $f\in F_{n}$ (recall that $\phi^{''}=\phi_{m_{1},s(n+1)}|_{P_{1}A_{m_{1}}P_{1}},\phi^{'}:=\phi_{m_{1},m_{2}}|_{P_{0}A_{m_{1}}P_{o}}$).

Combining with (4.2) and (4.4), and definition of $\beta_{n}$ and $\alpha_{n+1}$, the above inequality yields
$$\|\phi_{s(n),s(n+1)}(f)-(\alpha_{n+1}\circ\beta_{n+1})(f)\| < 9\varepsilon_{n}+3\varepsilon_{n}+2\varepsilon_{n}=14\varepsilon_{n},\;\forall\; f\in F_{n}.$$
Combining (1.3), the first inequality of (2.4) and the definition of $\beta_{n}$ and $\psi_{n,n+1}$, we have
$$\|\psi_{n,n+1}(g)-(\beta_{n}\circ\alpha_{n})(g)\|<10\varepsilon_{n}+3\varepsilon_{n}<14\varepsilon_{n},\;\forall\; g\in \widetilde{E}_{n}.$$
So we obtain(0.4).
The theorem follows from Proposition 4.1 of [GJLP2].\\
\end{proof}
Note that if $q\in M_{l}(I_{k})$, then $qM_{k}(I_{k})q$ isomorphic to $M_{l_{1}}(I_{k})$. Combining with the main theorem of [GJLP2] (see Theorem 4.2, and 2.7 of [GJLP2]) we have

\noindent\textbf{Theorem 3.2.} Suppose that $A=\lim(A_{n}=\oplus P_{n, i}M_{[n,i]}(C(X_{n,i}))P_{n,i})$ is an $AH$ inductive limit with $dim(X_{n,i})\leqslant M$ for a fixed positive integer $M$ such that limit algebra $A$ has ideal property. Then $A$ can be rewrite as inductive limit $\lim(B_{n}=\oplus B^{i}_{n}, \psi_{n,m})$, where either $B^{i}_{n}=Q_{n,i}M_{[n,i]^{'}}(C(Y_{n,i}))Q_{n,i}$ with $Y_{n,i}$ being one of the spaces $\{pt\},[0,1],S^{1},\{T_{\uppercase\expandafter{\romannumeral2},k}\}^{\infty}_{k=2}$ or $B^{i}_{n}=M_{[n,i]^{'}}(I_{l_{(n,i)}})$ a dimension drop algebra.

\clearpage

\begin{tiny}

\begin{small}

\end{small}

\end{tiny}

\end{document}